\documentclass[10pt,reqno]{amsart}
\usepackage{amsmath,amscd,amsfonts,amsthm,amssymb,latexsym,mathrsfs}
\usepackage{graphicx,epsf,epic,eepic,epsfig}
\usepackage{url}

\theoremstyle{plain}
   \newtheorem{theorem}{Theorem}[section]
   
   \newtheorem{lemma}[theorem]{Lemma}
   \newtheorem{corollary}[theorem]{Corollary}
   \newtheorem{problem}[theorem]{Problem}

\theoremstyle{definition}
   
   \newtheorem{example}[theorem]{Example}

\theoremstyle{remark}
   \newtheorem{remark}[theorem]{Remark}

\author[P.~Br\"and\'en]{Petter Br\"and\'en}
\address{Department of Mathematics, Royal Institute of Technology, SE-100 44 Stockholm,
Sweden}
\email{pbranden@math.kth.se}



\thanks{PB is a Royal Swedish Academy of Sciences Research Fellow
  supported by a grant from the Knut and Alice Wallenberg
  Foundation.}
\def\kkk{\kern.2ex\mbox{\raise.5ex\hbox{{\rule{.35em}{.12ex}}}}\kern.2ex}

\numberwithin{equation}{section}

\newcommand{\kk}{\mathbf{k}}
\newcommand{\LP}{\mathscr{L{\kkk}P}}
\newcommand{\LY}{\mathscr{L{\kkk}Y}}
\newcommand{\NN}{\mathbb{N}}

\newcommand{\PP}{\mathscr{P}}

\newcommand{\FF}{\mathscr{F}}

\newcommand{\RR}{\mathbb{R}}
\newcommand{\CC}{\mathbb{C}}

\newcommand{\KK}{\mathbb{K}}

\renewcommand{\Im}{{\rm Im}}
\renewcommand{\Re}{{\rm Re}}

\def\newop#1{\expandafter\def\csname #1\endcsname{\mathop{\rm
#1}\nolimits}}

\newop{per}
\newop{roots}
\newop{edges}
\newop{diag}
\newop{supp}
\newop{Proj}
\newop{JP}
\newop{rank}
\newop{Sym}
\newop{sign}
\newop{Int}
\newop{Cap}
\newop{PolO}
\newop{St}
\newop{glb}
\newop{MAP}
\newop{MOD}

\begin{document}

\title[Lee--Yang and P\'olya--Schur programs. III.]{The Lee--Yang and P\'olya--Schur programs. III. \\ Zero-preservers on Bargmann--Fock spaces}

\maketitle

\thispagestyle{empty}
\begin{abstract}
We characterize linear operators preserving zero-restrictions on entire functions in weighted Bargmann--Fock spaces.  The characterization  extends the main results of \cite{BB1,BBI,BBII} to the realm of entire functions, and translates into an optimal, albeit formal, Lee--Yang theorem. 
\end{abstract}

\section{Introduction}
The problem of describing linear operators preserving zero-restrictions on polynomials and transcendental entire functions has since the pioneering work of Hermite, Laguerre, Jensen and P\'olya been revisited frequently, see e.g. \cite{BB1,BBI, CC1, Le, LS}.
Already in the seminal papers \cite{LYII, LYI} of Lee and Yang it was made evident that linear operators preserving zero restrictions play a prominent role in understanding phase transitions of spin models in statistical physics. The method of using such preservers  was further developed by Lieb and Sokal \cite{LS} to prove a general Lee--Yang theorem.  In  a series of papers \cite{BB1,BBI,BBII} joint with Borcea we have characterized linear operators on polynomials preserving the property of being non-vanishing whenever the variables are in prescribed open circular regions. This constitutes a vast generalization of P\'olya and Schur's theorem \cite{PS} characterizing diagonal linear operators on polynomials preserving real--rootedness.   In this paper we extend the main results of \cite{BB1,BBI, BBII} to weighted Bargmann--Fock spaces of entire functions. The extension makes the connection between the P\'olya--Schur and Lee--Yang programs truly transparent.  Indeed, our characterization (Theorem \ref{LP-ext}) of Laguerre--P\'olya preservers translates directly into an optimal, albeit formal, Lee--Yang theorem (Theorem \ref{snurr}). 

\section{Laguerre--P\'olya preservers}

We say that a polynomial $P(z) \in \CC[z_1,\ldots, z_n]$ is \emph{stable} if $P(z) \neq 0$ whenever $z \in H^n$ where 
$H=\{ z \in \CC : \Im(z) >0\}$, and that an entire function $f(z)$ in $n$ variables is in the \emph{complex Laguerre--P\'olya class}, $f(z) \in \LP_n(\CC)$, if $f$ is the uniform limit on compact subsets of $\CC^n$ of stable polynomials. The (real) \emph{Laguerre--P\'olya class}, 
$\LP_n(\RR)$, consists of those functions in $\LP_n(\CC)$ with real coefficients. Laguerre and P\'olya proved that a univariate entire function is 
in the Laguerre--P\'olya class if and only if it may be expressed as 
\begin{equation}\label{lpform}
f(z)= Cz^ne^{az-bz^2}\prod_{k=1}^\omega (1+ x_kz)e^{-x_kz}, 
\end{equation}
where $C,a, x_k \in \RR$ for all $k$, and $b \geq 0$, $\omega \in \NN\cup\{\infty\}$, $n \in \NN$ and $\sum_k x_k^2 < \infty$.

The symbol of a linear operator $T : \CC[z_1,\ldots, z_n] \rightarrow \CC[z_1,\ldots,z_m]$ is the formal power series
$$
G_T(z,w)= T(e^{z\cdot w}):= \sum_{\alpha \in \NN^n} T(z^\alpha)\frac {w^\alpha} {\alpha!},
$$
where $\alpha! = \alpha_1! \cdots \alpha_n!$, $z^\alpha = \prod_{i=1}^n z_i^{\alpha_i}$ and $z\cdot w= z_1w_1+\cdots+z_nw_n$. We say that $T$ \emph{preserves stability} if 
$T(P)$ is stable or identically zero whenever $P$ is stable. The following characterizations of stability preservers was achieved in \cite{BBI}. 

\begin{theorem}\label{trans-C}
Let $T : \CC[z_1,\ldots, z_n] \rightarrow \CC[z_1,\ldots,z_m]$ be a linear operator. Then $T$ preserves stability if and only if 
\begin{enumerate}
\item The rank of $T$ is at most one and $T$ is of the form 
$$
T(P)= \alpha(P)Q, 
$$
where $\alpha : \CC[z_1,\ldots, z_n] \rightarrow \CC$ is a linear functional and $Q$ is a stable polynomial, or 
\item $G_T(z,-w) \in \LP_{m+n}(\CC)$, where $-w=(-w_1,\ldots, -w_n)$. 
\end{enumerate}
\end{theorem}

\begin{theorem}\label{trans-R}
Let $T : \RR[z_1,\ldots, z_n] \rightarrow \RR[z_1,\ldots,z_m]$ be a linear operator. Then $T$ preserves real stability if and only if 
\begin{enumerate}
\item The rank of $T$ is at most two and $T$ is of the form 
$$
T(P)= \alpha(P)Q + \beta(P)R, 
$$
where $\alpha, \beta : \RR[z_1,\ldots, z_n] \rightarrow \RR$ are linear functionals and $Q+iR$ is a stable polynomial, or 
\item $G_T(z,w) \in \LP_{m+n}(\RR)$, or
\item $G_T(z,-w) \in \LP_{m+n}(\RR)$. 
\end{enumerate}
\end{theorem}
Since a real univariate polynomial is stable if and only if it has only real zeros, Theorem \ref{trans-R} characterizes real zero preservers when $n=m=1$.

We want to extend linear stability preservers to act on entire functions. More precisely we want entire functions in the (complex) Laguerre--P\'olya class to be mapped on entire functions in the (complex) Laguerre--P\'olya class. To achieve this we should at least demand that stable polynomials should be mapped into $\LP_m(\CC)$. However we shall see that this weakest requirement still allows us to extend the domain to classes of entire functions of bounded growth.  Let $\KK[[z_1,\ldots, z_n]]$, where $\KK=\RR$ or $\CC$, be the space of all formal power series with coefficients in $\KK$.  A $\KK$-linear operator $T : \KK[z_1,\ldots, z_n] \rightarrow  \KK[[z_1,\ldots, z_m]]$ is a called \emph{Laguerre--P\'olya preserver}
if 
$$
T(\LP_n(\KK) \cap \KK[z_1,\ldots, z_n]) \subseteq \LP_{m}(\KK), 
$$
that is, if stable polynomials are mapped into the (complex) Laguerre--P\'olya class. 
The symbol of a linear operator $T : \KK[z_1,\ldots, z_n] \rightarrow \KK[[z_1,\ldots,z_m]]$ is the formal power series
$$
G_T(z,w)= T(e^{z\cdot w}):= \sum_{\alpha \in \NN^n} T(z^\alpha) \frac {w^\alpha} {\alpha!}.
$$
Theorems \ref{trans-C} and \ref{trans-C} extend naturally to this general setting. 
\begin{theorem}\label{LP-C}
Let $T : \CC[z_1,\ldots, z_n] \rightarrow \CC[[z_1,\ldots,z_m]]$ be a linear operator. Then $T$ is a Laguerre--P\'olya preserver if and only if 
\begin{enumerate}
\item The rank of $T$ is at most one and $T$ is of the form 
$$
T(P)= \alpha(P)f, 
$$
where $\alpha : \CC[z_1,\ldots, z_n] \rightarrow \CC$ is a linear functional and $f \in \LP_m(\CC)$, or 
\item $G_T(z,-w) \in \LP_{m+n}(\CC)$.
\end{enumerate}
\end{theorem}
\begin{lemma}\label{somev}
Let $f(z,w)$ be a formal power series in $z_1,\ldots, z_m, w_1,\ldots, w_n$. Write $f$ as 
$$
f(z,w)=\sum_{\alpha \in \NN^n} a_\alpha(z) w^\alpha.
$$
Then $f \in \LP_{m+n}(\CC)$ if and only if for each $\beta \in \NN^n$ 
$$
\Lambda_\beta(f):=\sum_{\alpha \in \NN^n} (\beta)_\alpha a_\alpha(z) w^\alpha \in \LP_{m+n}(\CC),
$$
where 
$$
(\beta)_\alpha := \alpha! \prod_{i=1}^m \binom {\beta_i}{\alpha_i}.
$$
\end{lemma}

\begin{proof}
The lemma was proved in \cite[Theorem 6.1]{BBI} in the case when $a_\alpha$ is a polynomial for all $\alpha$. However the lemma follows from the special case. By \cite[Theorem 6.1]{BBI} 
\begin{align*}
f(z,w) \in \LP_{m+n}(\CC) &\Longleftrightarrow \Lambda_{\gamma \oplus \beta}(f) \in \LP_{m+n}(\CC) \mbox{ for all } \gamma \in \NN^m, \beta \in \NN^n \\
&\Longleftrightarrow \Lambda_{\beta}(f) \in \LP_{m+n}(\CC) \mbox{ for all }  \beta \in \NN^n. 
\end{align*}

\end{proof}

\begin{proof}[Proof of Theorem \ref{LP-C}]
For $\beta \in \NN^m$ let $T_\beta : \CC[z_1,\ldots, z_n] \rightarrow \CC[z_1,\ldots,z_m]$ be defined by $T_\beta=\Lambda_\beta \circ T$.
 Clearly $T$ has rank at most one if and only if $T_\beta$ has rank at most one for all $\beta$ with $\min_{1\leq i \leq n} \beta_i$ sufficiently large. By Lemma \ref{somev} $T$ is a Laguerre--P\'olya preserver if and only if $T_\beta$ preserves stability for all $\beta \in \NN^m$. Since $G_{T_\beta}(z,w)= \Lambda_\beta(G_T)$, where $\Lambda_\beta$ acts on the $z$-variables, the theorem follows from Theorem \ref{trans-C} and Lemma \ref{somev}. 
\end{proof}
The proof of the real version of Theorem \ref{LP-C} follows similarly, and is therefore omitted.
\begin{theorem}\label{LP-R}
Let $T : \RR[z_1,\ldots, z_n] \rightarrow \RR[[z_1,\ldots,z_m]]$ be a linear operator. Then $T$ is a Laguerre--P\'olya preserver if and only if 
\begin{enumerate}
\item The rank of $T$ is at most two and $T$ is of the form 
$$
T(P)= \alpha(P)f + \beta(P)g, 
$$
where $\alpha, \beta : \RR[z_1,\ldots, z_n] \rightarrow \RR$ are linear functionals and $f+ig \in \LP_m(\CC)$, or 
\item $G_T(z,w) \in \LP_{m+n}(\RR)$, or
\item $G_T(z,-w) \in \LP_{m+n}(\RR)$. 
\end{enumerate}
\end{theorem}

\begin{remark}\label{formalad}
If $T : \KK[z_1,\ldots, z_n] \rightarrow \KK[[z_1,\ldots, z_m]]$ is a linear operator we may define the \emph{formal adjoint},  
$T^\# : \KK[z_1,\ldots, z_m] \rightarrow \KK[[z_1,\ldots, z_n]]$,   as the linear operator with symbol $G_{T^\#}(z,w) = \overline{G_T(\overline{w},\overline{z})}$. Since $(v,w) \in H \times (-H)$ if and only if $(\overline{w}, \overline{v}) \in H \times (-H)$ we see that $T$ is a Laguerre--P\'olya preserver if and only if $T^\#$ is a Laguerre--P\'olya preserver (provided that $T$ is of rank greater than one if $\KK=\CC$, and greater than two if $\KK=\RR$). This duality can be seen as a vast generalization of a famous theorem due to Hermite, Jensen, P\'olya and Poulain:  
Let 
$T$ be a formal differential operator with constant coefficients of the form $T=g(d/dz)$ where $g$ is a real formal power series. Then $T$ preserves real--rootedness if and only if $g \in \LP_1(\RR)$. The formal adjoint of $T$ is the operator defined by $T^\#(f)= g(z)f(z)$, and so $T^\#$ is a Laguerre--P\'olya preserver if and only if $g \in \LP_1(\RR)$. 

In the next section we shall see that the formal adjoints considered here are actually proper adjoints in Hilbert spaces. 
\end{remark}

\section{Preservers on weighted Bargmann--Fock spaces}
A priori, the linear operators in Theorems \ref{LP-C} and \ref{LP-R} may only be applied to polynomials. However we shall see that the domain extends naturally to entire functions of bounded growth.  We want to find conditions on $G_T(z,w)$ that allow us to extend the domain to spaces 
of entire functions. It follows from  \cite[Theorem 6.6]{BBI} that for each $f \in \LP_n(\CC)$ there are constants $A, B >0$ such that 
\begin{equation}\label{r-bound}
|f(z)| \leq Ae^{Br^2}, \quad \mbox{ whenever } |z_j| \leq r \mbox{ for all } 1\leq j \leq r.
\end{equation}
Hence functions in the Laguerre--P\'olya class are of order at most two and of bounded type. Lieb and Sokal \cite{LS} worked with certain Frech\'et spaces of entire functions, we find it more convenient to work with Hilbert spaces:  
For  $\beta \in \RR_+^n:= (0,\infty)^n$, define the $\beta$-\emph{weighted Bargmann--Fock space}\footnote{Bargmann--Fock spaces have many names, usually a combination of Bargmann, Fischer, Fock and Segal.}, $\FF_\beta$, to be the space of all entire functions $f(z)= \sum_{\alpha \in \NN^n}a_\alpha z^\alpha$ such that 
$$
\|f\|_\beta^2 = \sum_{\alpha \in \NN^n}\frac {\alpha!}{\beta^{\alpha}} |a_\alpha|^2 < \infty,  
$$
see \cite{Ba}. One may also write 
$$
\|f\|_\beta^2 = \int_{\CC^n} |f(z)|^2d\sigma_\beta(z):= \frac {\beta_1\cdots \beta_n}{\pi^{n}}\int_{\CC^n} |f(z)|^2 \exp\left(-\sum_{i=1}^n\beta_i|z_i|^2\right)dm
$$
where $m$ is Lebesgue measure on $\CC^n=\RR^{2n}$. By \eqref{r-bound} we see that each $f \in \LP_n(\CC)$ is in $\FF_\beta$ for some $\beta \in \RR_+^n$. To be more precise if the entire function $f$ is $\mathcal{O}(\exp\left(\beta_1 |z_1|^2/2  + \cdots + \beta_n |z_n|^2/2\right))$, then $f \in \FF_\gamma$ for all $\gamma \gg \beta$ (by which we mean $\gamma_j> \beta_j$ for all $1\leq j \leq n$). The space $\FF_\beta$ is a Hilbert space with inner product given by 
$$
\langle f,g \rangle_\beta = \sum_{\alpha \in \NN^n} \frac {\alpha!} {\beta^\alpha}a_\alpha \overline{b_\alpha}= \int_{\CC^n} f(z)\overline{g(z)} d\sigma_\beta(z), 
$$
and orthonormal basis
$$\left\{\sqrt{\frac {\beta^\alpha} {\alpha!}}z^\alpha \right\}_{\alpha \in \NN^n}.
$$
It has a \emph{reproducing kernel} given by $e_\beta(z,\overline{w}):= \exp\left(\sum_{j=1}^n \beta_j z_j \overline{w_j}\right)$, that is, 
\begin{equation}\label{repr}
f(w)= \langle f(z), e_\beta(z,\overline{w}) \rangle_\beta,  
\end{equation}
for all $f \in \FF_\beta$. In particular, by Cauchy--Schwartz inequality, 
\begin{align}\label{tou}
&|f(w)|^2 = |\langle f(z), e_\beta(z,\overline{w}) \rangle_\beta|^2 \leq \|f\|_\beta^2 \|e_\beta(z,\overline{w})\|_\beta^2 \\
& \leq \|f\|_\beta^2 \frac {\beta_1\cdots \beta_n}{\pi^{n}}\int_{\CC^n} \exp\left(-\sum_{i=1}^n\beta_i(|z_i|^2-2|z_i||w_i|)\right)dm = C(|w_1|, \ldots, |w_n|)\|f\|_\beta^2, \nonumber
\end{align}
and hence convergence in $\| \cdot \|_\beta$ implies uniform convergence on compact subsets of $\CC^n$. If $\alpha\in \RR_+^m$ and  $\beta \in \RR_+^n$, let $\alpha \oplus \beta =(\alpha_1,\ldots, \alpha_m,\beta_1,\ldots, \beta_n)$, $\alpha^{-1}=(\alpha_1^{-1}, \ldots, \alpha_m^{-1})$, and if $n=m$ let 
$\alpha\beta=(\alpha_1\beta_1, \ldots, \alpha_n\beta_n)$.  

Let $\LP_\beta(\CC)= \LP_n(\CC) \cap \FF_\beta$. The following theorem tells us to which weighted Bargmann--Fock spaces a Laguerre--P\'olya preserver may be extended. Theorem \ref{LP-ext} (1) is sharp and the converse is given by the last sentence of Theorem \ref{symbol-boundop}. 

\begin{theorem}\label{LP-ext}
Let $T : \CC[z_1, \ldots, z_n] \rightarrow \CC[[z_1,\ldots, z_m]]$ be a linear operator of rank at least two. Then $T$ is a Laguerre--P\'olya preserver if and only if $G_T(z,-w) \in \LP_{\beta \oplus \gamma}(\CC)$ for some $\beta \in \RR_+^m$ and $\gamma \in \RR_+^n$. 

Moreover if $G_T(z,-w) \in \LP_{\beta \oplus \gamma}(\CC)$, then 
\begin{enumerate}
\item $T$ extends to a bounded 
linear operator $T : \FF_{\alpha} \rightarrow \FF_\beta$ of the form \eqref{g-bound} and \eqref{g-int} for all $\alpha \leq \gamma^{-1}$, and  
\item $T : \LP_\alpha(\CC) \rightarrow \LP_\beta(\CC)$, for all $\alpha \leq \gamma^{-1}$.  
\end{enumerate}
\end{theorem}
           
The real version of Theorem \ref{LP-ext} is similar. In (1) below one may either consider $T$ as the obvious complexification of $T$, or consider $T$ as a 
linear operator on the real weighted Bargmann--Fock space. 

\begin{theorem}\label{LPR-ext}
Let $T : \RR[z_1, \ldots, z_n] \rightarrow \RR[[z_1,\ldots, z_m]]$ be a linear operator of rank at least three. Then $T$ is a Laguerre--P\'olya preserver if and only if $G_T(z,w) \in \LP_{\beta \oplus \gamma}(\RR)$ or $G_T(z,-w) \in \LP_{\beta \oplus \gamma}(\RR)$ for some $\beta \in \RR_+^m$ and $\gamma \in \RR_+^n$. 

Moreover if $G_T(z,\pm w) \in \LP_{\beta \oplus \gamma}(\RR)$, then 
\begin{enumerate}
\item $T$ extends to a bounded 
linear operator $T : \FF_{\alpha} \rightarrow \FF_\beta$ of the form \eqref{g-bound} and \eqref{g-int} for all $\alpha \leq \gamma^{-1}$, and  
\item $T : \LP_\alpha(\RR) \rightarrow \LP_\beta(\RR)$, for all $\alpha \leq \gamma^{-1}$.  
\end{enumerate}
\end{theorem}

\begin{example}\label{lp1}
 For $n=1$ and $\KK=\RR$, Laguerre--P\'olya preservers are linear operators that send polynomials with only real zeros into the Laguerre--P\'olya class. If the symbol of 
$T$ is in $\FF_{\beta \oplus \gamma}$, then $T$ extends to all $f$ as in \eqref{lpform} with $2b< 1/\gamma$. If $T$ is a multiplier sequence, see 
\cite{CC1, PS}, then its symbol is of the form 
$$
G_T(z,w)= C z^nw^n e^{\pm azw}\prod_{j=1}^\omega (1\pm x_j zw), 
$$
where $C \in \RR$, $a\geq 0$, $n \in \NN$, $\omega \in \NN\cup\{\infty\}$, $x_j > 0$ for all $j \in \NN$ and $\sum_{j}x_j < \infty$. Since 
$$
\prod_{j=1}^\omega (1+ x_j |z||w|) \leq \exp\left(|z||w|\sum_{j=N}^\infty x_j\right) \prod_{j=1}^{N-1}(1+ x_j |z||w|)
$$
for all $N \in \NN$ and 
$$
  2a|z||w| \leq \frac 1 s |w|^2 + a^2s|z|^2
$$
for all $s >0$, we see that $G_T(z,w) \in \FF_{(\beta, \gamma)}$ for all $\beta >a^2s$ and $\gamma > 1/s$. Hence $T : \FF_c \rightarrow \FF_d$ whenever $c >0$ and $d> a^2c$. 
\end{example}

\begin{theorem}\label{symbol-boundop}
Let $T : \CC[z_1,\ldots, z_n] \rightarrow \CC[[z_1,\ldots, z_m]]$ be a linear operator such that $G_T(z,w) \in \FF_{\beta \oplus \gamma}$. Then 
$T$ defines a bounded operator $T : \FF_\alpha \rightarrow \FF_\beta$ for all $\alpha \leq \gamma^{-1}$:
\begin{equation}\label{g-bound}
\|T(f)\|_\beta \leq \|G_T(z, \alpha w)\|_{\beta \oplus \alpha} \|f\|_\alpha.
\end{equation}
Moreover $T$ has the integral representation 
\begin{equation}\label{g-int}
T(f)(z)= \int_{\CC^n} f(w)G_T(z, \alpha\overline{w})d\sigma_\alpha(w).
\end{equation}
Conversely if $T : \FF_\alpha \rightarrow \FF_\beta$ is a bounded operator, then $G_T(z,w) \in \FF_{\beta \oplus \gamma}$ for all $\gamma \gg \alpha^{-1}$. 
\end{theorem}
\begin{proof} 
Suppose that $H(z,w) \in \FF_{\beta \oplus \alpha}$, and define a linear operator by 
$$
T(f)(z)= \int_{\CC^n} f(w)H(z, \overline{w})d\sigma_\alpha(w).
$$
This is well defined since $H(z_0,w) \in \FF_\alpha$ for each $z_0 \in \CC^n$, and then  by Cauchy--Schwartz inequality 
\begin{equation}\label{TH}
|T(f)(z)| \leq \|f\|_\alpha \|H(z,\cdot)\|_\alpha. 
\end{equation}
It follows from \eqref{repr} that 
$$
H(z,w)= T\left(e_\alpha(z,w))\right)= \sum_{\eta \in \NN^n} T(z^\eta)\frac {w^\eta \alpha^\eta}{\eta!} = G_T(z, \alpha w).
$$
By \eqref{TH} 
$$
\|T(f)\|_\beta^2 \leq \|f\|_\alpha^2 \int_{\CC^m} \|H(z,\cdot)\|_\alpha^2 d\sigma_\beta(z)= \|f\|_\alpha^2 \|H\|_{\beta\oplus \alpha}^2 < \infty. 
$$
Hence $T : \FF_\alpha \rightarrow \FF_\beta$ is a bounded operator. 

To prove the first part of the theorem we want to determine when $G_T(z,\alpha w) \in \FF_{\beta \oplus \alpha}$ given that 
$G_T(z,w) \in \FF_{\beta \oplus \gamma}$. However $G_T(z, \alpha w) \in \FF_{\beta \oplus \alpha^2\gamma}$ which implies 
$G_T(z,\alpha w) \in \FF_{\beta \oplus \alpha}$ whenever $\alpha^2 \gamma \leq \alpha$ from which the first part of the theorem follows.

Conversely suppose that $T : \FF_\alpha \rightarrow \FF_\beta$ is a bounded operator so that 
$
\|T(f)\|_\beta \leq C \|f\|_\alpha
$
for all $f \in \FF_\alpha$ and some $C>0$. For each fixed $w \in \CC^n$ define an entire function in $z$ by  $H(z,w)= T(e_{\alpha}(z,w))$. Now $H(z,w)$ defines an entire function on $\CC^{n+m}$ as can seen as follows. Let $E_k(z,w)= \prod_{j=1}^n(1+\alpha_jz_jw_j/k)^k$, then 
\begin{align*}
&\|T(E_k(z,w))\|_\beta^2 \leq C^2 \|E_k(z,w)\|_\alpha^2 \leq C^2 \|E_k(|z|,|w|)\|_\alpha^2 \\
&\leq D \int_{\CC^n} \exp\left( -\sum_{j=1}^n\alpha_j(|z_j|^2-2|z_j||w_j|)\right)dm = K(|w_1|, \ldots, |w_n|)<\infty. 
\end{align*}
By \eqref{tou} $T(E_k(z,w))$ is locally bounded, so by Vitali's theorem $T(E_k(z,w)) \rightarrow H(z,w)$ uniformly on compact subsets of $\CC^{n+m}$.  Hence $H(z,w)=T(e_\alpha(z,w))= G_T(z,\alpha w)$ is an entire function. 
Now let $\alpha' \gg \alpha$, then 
\begin{align*}
\|H(z,w)\|_{\beta\oplus\alpha'}^2 &= \int_{\CC^n} \|H(\cdot, w) \|_\beta^2 d\sigma_{\alpha'}(w)=\int_{\CC^n} \|T(e_\alpha(z, w)) \|_\beta^2d\sigma_{\alpha'}(w)\\ 
& \leq C^2 \int_{\CC^n} \| e_\alpha(z, w)\|_\alpha^2 d\sigma_{\alpha'}(w)= C^2 \|e_\alpha(z, w) \|_{\alpha \oplus \alpha'}^2 \\
&= C^2 \prod_{j=1}^n \frac {\alpha_j'}{\alpha_j'-\alpha_j}. 
\end{align*}
Hence $H(z,w) = G_T(z, \alpha w) \in \FF_{\beta \oplus \alpha'}$ for all $\alpha' \gg \alpha$, from which the theorem follows. 
\end{proof}

\begin{proof}[Proof of Theorem \ref{LP-ext}]
By Theorems \ref{trans-C} and \ref{symbol-boundop}, and \eqref{r-bound} it remain to prove (2). Let $f = \sum_{\gamma}a_\gamma z^\gamma \in \LP_\alpha(\CC)$, and let 
$$
f_k(z)= \sum_{\gamma \leq \kk} \frac {(\kk)_\gamma}{ \kk^\gamma} a_\gamma z^\gamma, 
$$
where $\kk = (k,\ldots, k)\in \NN^n$. By Lemma \ref{somev} $f_k$ is stable or identically zero for all $k$. Since ${(\kk)_\gamma}/{\kk^\gamma} \leq 1$ for all $k$ and $\gamma$ and 
${(\kk)_\gamma}/{\kk^\gamma} \rightarrow 1$ as $k \rightarrow \infty$ for all $\gamma\in \NN^n$ we have  $f_k \rightarrow f$ in $\FF_\alpha$. Since $T$ is bounded 
$T(f_k) \rightarrow T(f)$ in $\FF_\beta$, and hence $T(f_k) \rightarrow T(f)$ uniformly on compact subsets of $\CC^{n}$. Now $T(f_k) \in \LP_\beta(\CC)$  for all $k$ since $T$ is a Laguerre--P\'olya preserver, and the theorem follows. 
\end{proof}

The operators in question are bounded, so the dual operator is well-defined and bounded. Since  $T : \FF_\alpha \rightarrow \FF_\beta$ and $T^*:  \FF_\beta \rightarrow \FF_\alpha$ are related by 
$
\langle T(f),g \rangle_\beta =\langle f,T^*(g) \rangle_\alpha, 
$
we see that 
\begin{equation}\label{dualsymb}
G_{T^*}(w, \beta v)= \overline{G_{T}(\overline{v}, \overline{\alpha w})}, 
\end{equation}
by setting $f(z)=e_\alpha(z, \overline{w})$ and $g(z)=e_\beta(z, \overline{v})$ and using \eqref{repr}. As in Remark \ref{formalad} we obtain:
\begin{corollary}\label{duality}
Let $T : \FF_\alpha \rightarrow \FF_\beta$ be a bounded linear operator of rank at least two (or at least three if $T$ is considered as acting on real entire functions). Then $T$ is a  Laguerre--P\'olya preserver if and only if its 
dual $T^* :  \FF_\beta \rightarrow \FF_\alpha$ is a Laguerre--P\'olya preserver. 
\end{corollary}

\begin{example}\label{diffop}
Let $T$ be a differential operator with constant coefficient, $T= g(\partial / \partial z)$, where $\partial / \partial z=(\partial / \partial z_1, \ldots, \partial / \partial z_n)$ and $g$ is an entire function. For which $\gamma, \beta \in \RR_+^n$ is $T : \FF_\alpha \rightarrow \FF_\gamma$ a bounded operator? This was answered in \cite[Proposition 2.5]{LS}, and we shall see how it may be derived from Theorem \ref{symbol-boundop}.

The symbol of $T$ is $g(w)e^{z\cdot w}$, and $\|e^{z\cdot w}\|_{\beta \oplus \eta}^2= \prod_{j=1}^n(\beta_j\eta_j/(\beta_j\eta_j-1))<\infty$ if and only if $\beta \gg 1/\eta$. Hence if $\alpha \in \RR_+^n$ and 
$$
M_\alpha(g) :=  \sup_{z \in \CC^n}\left[ \exp\left(-\sum_{j=1}^n \alpha_j |z_j|^2/2\right)|g(z)|\right] < \infty, 
$$
then 
\begin{align*}
\|g(w)e^{z\cdot w}\|_{\beta \oplus (\eta+\alpha)}^2 &= C^2 \int_{\CC^{2n}} |g(w) e_{\alpha/2}(-w, \overline{w})|^2|e^{z\cdot w}|^2 d\sigma_{\beta \oplus \eta}(z,w)\\
&\leq C^2M_\alpha(g)^2 \|e^{z\cdot w}\|_{\beta \oplus \eta}^2, 
\end{align*}
where $C>0$ is a constant. By Theorem \ref{symbol-boundop}, $T : \FF_{1/(\eta+\alpha)} \rightarrow \FF_\beta$ whenever 
$\beta \gg 1/\eta$. Hence $T : \FF_\gamma \rightarrow \FF_\beta$ whenever $\gamma \ll 1/\alpha$ and $\beta \gg \gamma/(1-\alpha\gamma)$, where $1$ is the all ones vector. 

This is sharp by Theorem \ref{symbol-boundop}, which can also be seen from the following example of P\'olya: Let $g(z)=\exp(az/2)$ and $f(z)=\exp(bz/2)$ where $a,b>0$.  Then $M_\alpha(g)<\infty$ if and only if  $\alpha \geq a$ and $f \in \FF_\gamma$ if and only if $\gamma > b$. Hence if $ab<1$, then by the above,  
$g(d/dz) f \in \FF_\beta$ for any $\beta > b/(1-ab)$. Now 
$$
g(d/dz)f = \frac 1 {\sqrt{1-ab}}\exp\left(\frac 1 2 \frac b {1-ab} z^2\right), 
$$
so that $g(d/dz)f \in \FF_\beta$ if and only if $\beta > b/(1-ab)$.

\end{example}

\section{Lee--Yang theorems}
The Lee--Yang theorem and its extensions assert that the Fourier--Laplace transform of Gibbs measures of various spin models are nonzero whenever all variables are in the open right half-plane, and they serve as important tools in the rigorous study of phase transitions in lattice spin systems \cite{Ru1}.  We follow the approach to the Lee--Yang theorem developed by Lieb and Sokal \cite{LS} that uses linear operators preserving non-vanishing properties. For another successful method which uses Asano contractions we refer to \cite{Ru1,Ru2} and the references therein. 

Denote by $\LY_n(\CC)$ the space of all entire functions in $n$ variables that are 
uniform limits on compact subsets of $\CC^n$ of polynomials that are nonvanishing whenever all variables are in the open right half-plane of the complex plane. Thus $f(z_1, \ldots, z_n) \in \LY_n(\CC)$ if and only if $f(-iz_1, \ldots, -iz_n) \in \LP_n(\CC)$. 
A measure $\mu$ on $\RR^n$ is said to have the \emph{Lee--Yang property} if its Fourier--Laplace transform 
$$
\hat{\mu}(w) := \int_{\RR^n}e^{z\cdot w} d\mu(z)
$$
is an entire function in $\LY_n(\CC)$. More generally a continuous linear functional $\phi : \FF_\beta \rightarrow \CC$ has the Lee--Yang property if 
the map $\hat{\phi}$ defined by $w \mapsto \phi(e^{z\cdot w})$ defines an entire function in $\LY_n(\CC)$. It is natural to extend the definition to linear operators: A bounded linear operator $T : \FF_\alpha \rightarrow \FF_\beta$ has the Lee--Yang property if $T(e^{z\cdot w}) =G_T(z,w) \in \LY_{n+m}(\CC)$. 

\begin{example}\label{basex}
Here are a few basic examples of measures on $\RR$ with the Lee--Yang property: 
\begin{enumerate}
\item If $\mu= (\delta_a + \delta_b)/2$, where $\delta_a$ and $\delta_b$ are the Dirac measures centered at $a,b \in \RR$, then 
$$
\hat{\mu}(z) =\exp\left(\frac {a+b} 2 z \right) \cosh\left( \frac {a-b} 2 z \right). 
$$
Hence $\mu : \FF_c \rightarrow \CC$ has the Lee--Yang property for all $c>0$ and $a+b \geq 0$. Moreover if $a$ and $b$ are allowed to be non-real, then $\mu$ has the Lee--Yang property if and only if $a-b\in \RR$ and $\Re(a+b)\geq 0$. 
\item If $\mu$ is Lebesgue measure on the interval $[a,b]$, then 
$$
\hat{\mu}(z) = \int_a^b e^{zx}dx= \frac 2 z \exp\left(\frac {b-a} 2 z \right) \sinh\left( \frac {a+b} 2 z \right). 
$$
Hence $\mu : \FF_c \rightarrow \CC$ has the Lee--Yang property for all $c>0$. 
\item If $d\mu(x)= e^{-bx^2/2}dx$ on $\RR$ with $b>0$, then 
$$
\hat{\mu}(z)= \int_\RR e^{zx}e^{-bx^2/2} dx= \sqrt{2\pi/b}\exp\left(\frac {z^2}{2b}\right).
$$
Hence $\mu : \FF_a \rightarrow \CC$ has the Lee--Yang property for all $a < b$. 
\end{enumerate}
\end{example}

The results developed in the previous sections apply (by a change of variables) to \emph{Lee--Yang preservers} which we define to be linear operators 
$T : \CC[z_1,\ldots, z_n] \rightarrow \CC[[z_1,\ldots, z_m]]$ that map $\LY_n(\CC)\cap \CC[z_1,\ldots, z_n]$  
into $\LY_m(\CC)$. For $\beta \in \RR_+^n$, let $\LY_\beta(\CC)= \LY_n(\CC) \cap \FF_\beta$. 

\begin{theorem}\label{LY-ext}
Let $T : \CC[z_1, \ldots, z_n] \rightarrow \CC[[z_1,\ldots, z_m]]$ be a linear operator of rank at least two. Then $T$ is a Lee--Yang preserver if and only if $T$ has the Lee--Yang property, that is, $G_T(z,w) \in \LY_{\beta \oplus \gamma}(\CC)$ for some $\beta \in \RR_+^m$ and  $\gamma \in \RR_+^n$. 

Moreover if $G_T(z,w) \in \LY_{\beta \oplus \gamma}(\CC)$, then 
\begin{enumerate}
\item $T$ extends to a bounded 
linear operator $T : \FF_{\alpha} \rightarrow \FF_\beta$ of the form \eqref{g-bound} and \eqref{g-int} for all $\alpha \leq \gamma^{-1}$, and  
\item $T : \LY_\alpha(\CC) \rightarrow \LY_\beta(\CC)$, for all $\alpha \leq \gamma^{-1}$.  
\end{enumerate}
\end{theorem}

\begin{corollary}\label{LY-dual}
Let $T : \FF_\alpha \rightarrow \FF_\beta$ be a bounded linear operator  of rank at least two. Then $T$ is a Lee--Yang preserver if and only if its 
dual $T^* :  \FF_\beta \rightarrow \FF_\alpha$ is a Lee--Yang preserver. 
\end{corollary}

\begin{remark}\label{extendef}
If $T : \FF_\alpha \rightarrow \FF_\beta$ is a bounded linear operator, $z_1', \ldots, z_k'$ are new variables and $\gamma \in \RR_+^k$, then 
$T$ extends to a bounded linear operator $\tilde{T} : \FF_{\alpha \oplus \gamma} \rightarrow \FF_{\beta \oplus \gamma}$ by setting 
$\tilde{T}(f(z, z'))= T(f(z,z'))$ and where $T$ only acts on the $z$-variables. Note also that the symbol of $\tilde{T}$ is 
$e^{z' \cdot w'}T(e^{z \cdot w})$ so that $T$ has the Lee--Yang property if and only if $\tilde{T}$ has the Lee--Yang property. 
\end{remark}

The next theorem shows that bounded linear operators with the Lee--Yang property are closed under composition. This can be seen as an ultimate generalization of \cite[Proposition 2.9]{LS} and serves as a fundamental tool to prove Lee--Yang theorems for one-component models.

\begin{theorem}\label{snurr}
Suppose that $T : \FF_\alpha \rightarrow \FF_\beta$ and $S : \FF_\beta \rightarrow \FF_\gamma$ have the Lee--Yang property. Then so does 
$S \circ T : \FF_\alpha \rightarrow \FF_\beta$. 

In particular, if $\phi : \FF_\beta \rightarrow \CC$ and $T : \FF_\alpha \rightarrow \FF_\beta$ have the Lee--Yang property, then $\phi \circ T: \FF_\alpha \rightarrow \CC$ has the Lee--Yang property. 
\end{theorem}

\begin{proof}
The symbol of $S \circ T$ is $\tilde{S}(G_T(z,w))$, where $\tilde{S} : \FF_{\beta\oplus \kappa} \rightarrow \FF_{\gamma\oplus \kappa}$  is as in Remark \ref{extendef}. Now $G_T(z,w) \in \LY_{\beta\oplus \kappa}$ for all $\kappa \gg \alpha^{-1}$. By Theorem \ref{LY-ext} and Remark \ref{extendef} $\tilde{S}$ is a Lee--Yang preserver and thus $\widehat{S\circ T} = \tilde{S}(G_T(z,w)) \in \LY_{\gamma\oplus \kappa}$. 

\end{proof}
The following corollary is an equivalent formulation of the most general (formal) one component Lee--Yang theorem in \cite{LS} from which many others follow.   
\begin{corollary}\label{gls}
Suppose that $\phi : \FF_\beta \rightarrow \CC$ has the Lee--Yang property. If $\alpha, \beta, \gamma \in \RR_+^n$ satisfy $\alpha + \gamma \leq \beta$,  and $g \in \LY_n(\CC)$ satisfies 
$M_\alpha(g) < \infty$, then $\psi : \FF_\gamma \rightarrow \CC$ defined by 
$
\psi(f)= \phi(fg)
$
has the Lee--Yang property. 
\end{corollary}
\begin{proof}
Clearly the operator $T(f)=fg$ is a Lee--Yang preserver. By Theorem \ref{snurr}  it remains to prove that $T : \FF_\gamma \rightarrow \FF_{\alpha+\gamma}$ is a bounded operator. 
If $f \in \FF_\gamma$ and $M_\alpha(g) < \infty$, then 
\begin{align*}
\|gf\|_{\alpha+\gamma}^2 &= \prod_{j=1}^n(1+\alpha_j/\gamma_j) \int_{\CC^n} |g(z) e_{\alpha/2}(-z, \overline{z})|^2|f(z)|^2 d\sigma_{\gamma}(z)\\
&\leq \prod_{j=1}^n(1+\alpha_j/\gamma_j) M_\alpha(g)^2\|f\|_{\gamma}^2. 
\end{align*}
\end{proof}

Since $e^{z_1z_2}, e^{z_1^2} \in \LY_2(\CC)$, we see that $e_J(z)=\exp(\sum_{i,j=1}^nJ_{ij} z_iz_j) \in \LY_n(\CC)$ for all matrices $J$ such that $J_{ij} \geq 0$ for all $j$. The original Lee--Yang theorem \cite{LYII} states that measures of the type 
$$
d\mu = e_J(z)d\mu_1(z_1)\cdots d\mu_n(z_n)
$$
where $\mu_1, \ldots, \mu_n$ are measures on $\RR$ as in Example 4.1 (1) with $a=-b=1$ and $J$ is a (entry-wise) nonnegative symmetric matrix. Clearly the direct product of two measures with the Lee--Yang property has the Lee--Yang property, so the original Lee--Yang theorem follows from Corollary \ref{gls}.

Newman's Lee--Yang theorem \cite{New} asserts (in our language) that 
$$
d\mu = e_J(z)d\mu_1(z_1)\cdots d\mu_n(z_n)
$$
has the Lee--Yang property whenever $\mu_1, \ldots, \mu_n$ are even measures on $\RR$ with the Lee--Yang property and $J$ is a nonnegative symmetric matrix such that $\hat{\mu} \in \FF_\beta$ for some $\beta \in \RR_+^n$. Hence Newman's theorem also follows from Corollary \ref{gls}. 

The Lee--Yang theorem of Lieb and Sokal \cite[Theorem 3.2]{LS} asserts that 
$$
d\mu = e_J(z)d\mu_0
$$
has the Lee--Yang property whenever $\mu_0$ has the Lee--Yang property, $J$ is a nonnegative symmetric matrix, and $\hat{\mu} \in \FF_\beta$ for some 
$\beta \in \RR_+^n$. Hence this theorem also follows from Corollary \ref{gls}. 

\begin{remark}
To apply Newman's or Lieb and Sokal's theorem one needs to know for which symmetric matrices $A$ with nonnegative entries  
$\exp\left(\sum_{i,j}A_{ij}z_iz_j\right) \in \FF_\beta$, so that one can use Corollary \ref{gls}. This happens if and only if 
$$
\sum_{i,j}A_{ij}|z_i||z_j| \leq \sum_{i} \alpha_i |z_i|^2
$$
for some $\alpha \ll \beta/2$, that is, if and only if $\| D_\alpha^{-1/2} A D_\alpha^{-1/2}\| \leq 1$ for some $\alpha \ll \beta/2$, where $\|\cdot \|$ denotes the operator norm and $D_\alpha= \diag(\alpha_1,\ldots, \alpha_n)$ is a diagonal matrix. 
\end{remark}
\section{An open problem}

We end by generalizing an important open problem from \cite{LS}. Let $\Gamma \subset \RR^n$ be an open convex cone, and let 
$\PP_n(\Gamma)$ be the set of polynomials in $n$ variables that are non-vanishing whenever the real parts of the variables are in $\Gamma$. 
\begin{problem}\label{conepr}
Let $\Gamma \subset \RR^n$ and $\Lambda \subset \RR^m$ be be two open convex cones.  Characterize all linear operators $T : \CC[z_1,\ldots, z_n] \rightarrow \CC[z_1,\ldots, z_m]$ such that 
$T(\PP_n(\Gamma))  \subseteq \PP_m(\Lambda) \cup\{0\}$. 
\end{problem}

When $\Gamma = \RR_+^n$ and $\Lambda = \RR_+^m$, Problem \ref{conepr} is just Theorem \ref{trans-C} (by a rotation of the variables). 
A solution to Problem \ref{conepr} would entail optimal Lee--Yang theorems for $N$-component models when $N \geq 3$, see \cite[Section 5]{LS} where partial results on Problem \ref{conepr} for differential operators in the Weyl-algebra were obtained. Progress on Problem \ref{conepr} would also be interesting for the 
convex optimization community, see e.g. \cite{Re}, since a homogeneous polynomial $P$ is in $\PP_n(\Gamma)$ if and only if $P$ is hyperbolic with hyperbolicity cone containing $\Gamma$. Thus Problem \ref{conepr} asks how one may deform a hyperbolic polynomial and retain hyperbolicity.

\end{document}